\documentclass[12pt]{article}
\usepackage{amsfonts}

\usepackage{latexsym}
\usepackage{amssymb}
\usepackage{epsfig}
\usepackage{amsmath}
\usepackage{amsthm}
\usepackage[mathscr]{eucal}
\usepackage{subfigure}
\usepackage{graphicx}

\newtheorem{Lemma}{Lemma}

\newtheorem{Theorem}[Lemma]{Theorem}

\newtheorem{Corollary}[Lemma]{Corollary}

\renewcommand{\qed}{\hfill{\ \ \rule{2mm}{2mm}} \vspace{0.2in}}

\newcommand{\ind}{1\hspace{-2.3mm}{1}}

\begin{document}

\title{First Passage Percolation with nonidentical passage times}
\author{ \textbf{Ghurumuruhan Ganesan}
\thanks{E-Mail: \texttt{gganesan82@gmail.com} } \\
\ \\
EPFL, Lausanne }
\date{}
\maketitle

\begin{abstract}
In this paper we consider first passage percolation on the square lattice \(\mathbb{Z}^d\) with
passage times that are independent and have bounded \(p^{th}\) moment for some \(p > 6(1+d),\) but not necessarily identically distributed. For integer \(n \geq 1,\) let \(T(0,n)\) be the minimum time needed to reach the point \((n,\mathbf{0})\) from the origin. We prove that \(\frac{1}{n}\left(T(0,n) - \mathbb{E}T(0,n)\right)\) converges to zero in \(L^2\) and use a subsequence argument to obtain almost sure convergence. As a corollary, for i.i.d. passage times, we also obtain the usual almost sure convergence of \(\frac{T(0,n)}{n}\) to a constant \(\mu.\)

\vspace{0.1in} \noindent \textbf{Key words:} First passage percolation
nonidentical passage times.

\vspace{0.1in} \noindent \textbf{AMS 2000 Subject Classification:} Primary:
60J10, 60K35; Secondary: 60C05, 62E10, 90B15, 91D30.
\end{abstract}

\bigskip

\section{Introduction} \label{intro}
Consider the square lattice \(\mathbb{Z}^d\) with edges \(\{e_i\}_{i \geq 1}.\) The passage times \(\{t(e_i)\}_{i}\) are independent random variables that satisfy the following conditions.\\
(i) We have that \(\sup_{i} \mathbb{P}(t(e_i) < \epsilon)  \longrightarrow 0\) as \(\epsilon \downarrow 0.\)\\
(ii) There exists a constant \(\eta > 0\) such that \(\sup_{i} \mathbb{E}(t(e_i))^{6(1+d) + \eta} < \infty.\)\\

For \(n \geq 1,\) we are interested in the shortest time path from \((0,\mathbf{0})\) to \((n,\mathbf{0}),\) where \(\mathbf{0} \) is the \((d-1)-\)dimensional zero vector. To define such a path, we proceed as follows. For any fixed path \(\pi\) starting from the origin and containing \(k\) edges \(e_1,...,e_k,\) we define the passage time to be \(T(\pi) = \sum_{i=1}^{k} t(e_i).\) Using (ii), we get that there exists a constant \(0 < \beta_1 < \mu\) such that
\begin{equation}\label{t_pi}
\mathbb{P}(T(\pi) \leq \beta_1 k) \leq e^{-dk}
\end{equation}
for all \(k \geq 1.\) We prove all estimates at the end of this section. By (iii) we have that \(\mu = \sup_{i}\mathbb{E}t(e_i) <\infty\) and by (ii) we have that \(\mu \geq \inf_i \mathbb{E}t(e_i) > 0.\) Let \(E_{k}\) denote the event that there exists a path starting from \((0,\mathbf{0})\) containing \(r \geq \frac{8 \mu}{\beta_1}k\) edges and whose passage time is less than \(\beta_1 r.\) Since there are at most \((2d)^{r}\) paths containing \(r\) edges, we have that
\begin{equation}\label{a_0k}
\mathbb{P}(E_{k}) \leq \sum_{r \geq 8\mu\beta_1^{-1} k} (2d)^{r}e^{-dr}  \leq C e^{-\beta_2 k}
\end{equation}
for all \(k \geq 1\) and for some positive constants \(\beta_2\) and \(C.\) To obtain (\ref{a_0k}), we let \(\delta = d-\log(2d).\) Since \(de^{-d} \leq e^{-1} < \frac{1}{2}\) for all \(d \geq 2,\) we have that \(\delta > 0\) and we obtain 
\begin{eqnarray}
\mathbb{P}(E_k) \leq \sum_{r \geq 8\mu\beta_1^{-1} k} e^{-\delta r} = \frac{1}{1-e^{-\delta}} e^{-\delta 8\mu\beta_1^{-1} k}. \nonumber
\end{eqnarray} 

For \(i \geq 1,\) let \(f_i\) denote the edge between \((i-1,\mathbf{0})\) and \((i,\mathbf{0})\) and let \(A_n = \left\{\sum_{i=1}^{2n} t(f_i) \leq 6\mu n\right\},\) where \(\mu\) is as above. There exists a constant \(C_1 > 0\) such that
\begin{equation}\label{an_prob}
\mathbb{P}(A_n^c) \leq \frac{C_1}{n^2}
\end{equation}
for all \(n \geq 1.\)
Finally, setting \(F_n = E_{n}^c \cap A_n,\) we note that if \(F_n\) occurs, then the time taken to reach \((i,\mathbf{0})\) from \((0,\mathbf{0})\) is less than \(6\mu n,\) for each \(1 \leq i \leq 2n.\) Since \(E^c_{n}\) also occurs, every path starting from \((0,\mathbf{0})\) and containing \(r \geq \frac{8\mu}{\beta_1}n\) edges has passage time at least \(\beta_1 r \geq 8\mu n.\) Therefore, if \(F_n\) occurs, the shortest time path from \((0,\mathbf{0})\) to \((i,\mathbf{0})\) is contained in \(B_{8\mu\beta_1^{-1}n} := [-8\mu\beta_1^{-1}n,8\mu\beta_1^{-1}n]^d\) for each \(1 \leq i \leq n.\)

From (\ref{a_0k}) and (\ref{an_prob}), we have that
\begin{equation} \label{fn_prob}
\mathbb{P}(F_n^c) \leq \frac{C_2}{n^2}
\end{equation}
for some constant \(C_2 > 0\) and thus by Borel-Cantelli lemma, we have that \(\mathbb{P}(\liminf_n F_n)= 1.\) Fix \(\omega \in \liminf_n F_n\) and for every \(n \geq 1,\) define \(T(0,n)(\omega)\) to be the shortest time taken for reaching \((n,\mathbf{0})\) from \((0,\mathbf{0}).\) If there is more than one path that attains the shortest time, we provide an iterative procedure at the end of this section to choose a unique path.

We are interested in studying the convergence of \(\frac{T(0,n)}{n}.\) We have the following result.
\begin{Theorem} \label{thm1} We have that
\begin{equation}
\frac{1}{n}\left(T(0,n) - \mathbb{E}T(0,n)\right) \longrightarrow 0\;\;\text{a.s. and in }L^2
\end{equation}
as \(n \rightarrow \infty.\)
\end{Theorem}
For the case of independent and identically distributed (i.i.d.) random variables, we have the following Corollary.
\begin{Corollary} \label{cor1} If the passage times are i.i.d., we have that
\begin{equation}
\frac{T(0,n)}{n}  \longrightarrow \mu\;\;\text{a.s. and in }L^2
\end{equation}
as \(n \rightarrow \infty,\) for some constant \(\mu > 0.\)
\end{Corollary}
The constant \(\mu\) is also called the time constant; Alexander (1993), Cox and Durrett (1981), Kesten (1993) and Smythe and Wierman (2008) and references therein contain further material on first passage percolation.

The paper is organized as follows: In the rest of this section, we prove estimates (\ref{t_pi}) and (\ref{an_prob}) and provide an iterative procedure for choosing the minimum time path. In Section~\ref{pf1}, we prove Theorem~\ref{thm1} and  Corollary~\ref{cor1}.

To prove (\ref{t_pi}), we write \[\mathbb{P}(T(\pi) \leq \beta k)  = \mathbb{P}\left(\sum_{i=1}^{k}t(e_i) \leq \beta k\right)\] for a fixed \(\beta > 0.\) Since \(\{t(h_{i})\}_{i}\) are independent, we have for a fixed \(s  >0\) that
\begin{equation}\label{y_1_eq1}
\mathbb{P}(T(\pi) \leq \beta k) = \mathbb{P}\left(\sum_{i} t(e_i) \leq \beta k\right) \leq e^{s\beta k} \prod_{i=1}^{k}\mathbb{E}(e^{-st(e_i)}).
\end{equation}
For a fixed \(\epsilon > 0,\) we have that
\begin{eqnarray}
\mathbb{E}e^{-st(e_i)} &=& \int_{t(e_i) < \epsilon} e^{-st(e_i)} d\mathbb{P} + \int_{t(e_i) \geq \epsilon} e^{-st(e_i)} d\mathbb{P} \nonumber\\
&\leq& \int_{t(e_i) < \epsilon} e^{-st(e_i)} d\mathbb{P} + e^{-s\epsilon} \label{t_pi_an}\\
&\leq& \mathbb{P}(t(e_i) < \epsilon) + e^{-s\epsilon}. \nonumber
\end{eqnarray}
Using (i), the first term in the last expression is less than \(\frac{e^{-6d}}{2}\) if \(\epsilon > 0\) is small, independent of \(i.\) Fixing such an \(\epsilon,\) we choose \(s\) large so that the second term is also less than \(\frac{e^{-6d}}{2}.\) Substituting into (\ref{y_1_eq1}), we have that \[\mathbb{P}(T(\pi) \leq \beta k) \leq e^{s\beta k} e^{-3dk} \leq e^{-2dk},\] provided \(\beta > 0\) is small. We fix such a small \(\beta < \mu.\) 

To prove (\ref{an_prob}), we let \(\mu_i = \mathbb{E} t(f_i)\) and use Chebychev's inequality to write
\begin{equation}\label{a_0n_temp}
\mathbb{P}(A_n^c) \leq \mathbb{P}\left(\sum_{i=1}^{2n}X_i \geq 4\mu n\right) \leq \frac{1}{(4\mu n)^4}\mathbb{E}\left(\sum_{i} X_{i}\right)^{4},
\end{equation}
where \(X_i = t(f_i) - \mu_i.\) Since \(\{X_{i}\}_i\) are independent, we have that \(\mathbb{E}X_{i}X_{j} = 0\) for \(i \neq j.\) Thus we have \[\mathbb{E}\left(\sum_{i} X_{i}\right)^{4} = \sum_{i}\mathbb{E}X_{i}^4 + \sum_{i \neq j} \mathbb{E}X_{i}^2X_{j}^2 \leq C_1 n^2\] for some constant \(C_1 > 0\) by (ii). Substituting into (\ref{a_0n_temp}) proves (\ref{an_prob}).

Finally, we provide an iterative procedure to choose the shortest time path in the presence of multiple choices. For simplicity we provide for \(d  =2.\) An analogous procedure holds for general \(d.\) Fix \(\omega \in \liminf_n F_n\) and let \({\cal S}_1 = \{L_i\}_{1 \leq i \leq W} =\{(S_{i,1},...,S_{i,H_i})\}_{1 \leq i \leq W}\) be the set of all paths with the shortest passage time from \((0,0)\) to \((n,0).\) We note that \(W = W(\omega) < \infty.\) Let \(x_{i,j}\) and \(y_{i,j}\) be the \(x\)- and \(y\)-coordinates, respectively, of the centre of \(S_{i,j}.\) Let \(y'_1 = \min_{L_k \in {\cal S}_{1}} y_{k,1}\) and let \({\cal S}'_1 = \{L_k \in {\cal S}_1 : y_{k,1} = y'_1\}.\) Let \(x'_1 = \min_{L_k \in {\cal S}'_1} x_{k,1}.\) Let \(h_1\) be the edge attached to the origin whose centre has coordinates \((x'_1,y'_1).\) Clearly \(h_1\) is the first edge of some path in \({\cal S}'_1.\) Let \({\cal S}_2\) be the set of paths in \({\cal S}'_1\) whose first edge is \(h_1.\) Repeating the above procedure with \({\cal S}_2,\) we obtain an edge \(h_2\) attached to \(h_1.\) Continuing iteratively, this procedure terminates after a finite number of steps resulting in a unique path. Also, the final path obtained does not depend on the initial ordering of the paths.

\section{Proof of Theorem~\ref{thm1}}\label{pf1}
For \(n \geq 1,\) we define auxiliary random variables \(\{\hat{T}^{(n)}_k\}_{k \geq 1}\) defined as follows. For \(i \geq 1,\) let \(t_n(e_i) = \min(t(e_i), n^{\alpha}),\) where \(\alpha < \frac{1}{6}\) is a constant to be determined later. Since \(t_n(e_i) \leq t(e_i)\;\;\)a.s., we have that (i) and (ii) are satisfied by \(\{t_n(e_i)\}_i.\) For any fixed path \(\pi\) starting from the origin and containing \(k\) edges \(e_1,...,e_k,\) we define the passage time to be \(\hat{T}_n(\pi) = \sum_{i=1}^{k} t_n(e_i).\) We have
\begin{equation}\label{t_pi_p}
\mathbb{P}(\hat{T}_n(\pi) \leq \beta_1 k) \leq e^{-dk}
\end{equation}
for all \(k \geq 1.\) Here the constant \(\beta_1\) is the same as in (\ref{t_pi}) and is independent of \(n.\) To prove (\ref{t_pi_p}), we use the fact that \(\{t_n(e_{i})\}_{i}\) are independent and thus for a fixed \(s  >0\) we have that
\begin{equation}\label{y_1_eq1}
\mathbb{P}(\hat{T}_n(\pi) \leq \beta k) = \mathbb{P}\left(\sum_{i} t_n(e_i) \leq \beta k\right) \leq e^{s\beta k} \prod_{i=1}^{k}\mathbb{E}(e^{-st_n(e_i)}).
\end{equation}
For a fixed \(0 < \epsilon <1,\) we have that
\begin{eqnarray}
\mathbb{E}e^{-st_n(e_i)} &=& \int_{t_n(e_i) < \epsilon} e^{-st_n(e_i)} d\mathbb{P} + \int_{t_n(e_i) \geq \epsilon} e^{-st_n(e_i)} d\mathbb{P} \nonumber\\
&\leq& \int_{t_n(e_i) < \epsilon} e^{-st_n(e_i)} d\mathbb{P} + e^{-s\epsilon} \nonumber\\
&=& \int_{t(e_i) < \epsilon} e^{-st(e_i)} d\mathbb{P} + e^{-s\epsilon} \nonumber
\end{eqnarray}
which is the same as (\ref{t_pi_an}). The final equality is because \(\epsilon < 1\) and thus \(t_n(e_i)< \epsilon\) if and only if \(t(e_i) < \epsilon.\) Following an analogous analysis following (\ref{t_pi_an}) we obtain (\ref{t_pi_p}). For \(k \geq 1,\) let \(\hat{E}_{k}(n)\) denote the event that there exists a path \(\pi_1\) starting from \((0,\mathbf{0})\) containing \(r \geq \frac{8 \mu}{\beta_1}k\) edges and whose passage time \(\hat{T}_n(\pi_1)\) is less than \(\beta_1 r.\) As in (\ref{a_0k}) we have that
\begin{equation}\label{a_0k_p}
\mathbb{P}(\hat{E}_{k}(n)) \leq Ce^{-\beta_2 k}
\end{equation}
for all \(k \geq 1,\)  where \(\beta_2\) and \(C\) are as in (\ref{a_0k}).

As before, for \(i \geq 1\) let \(f_i\) denote the edge between \((i-1,\mathbf{0})\) and \((i,\mathbf{0})\) and for \(k \geq 1,\) let \(\hat{A}_n(k) = \left\{\sum_{i=1}^{2n} t_k(f_i) \leq 6\mu n\right\},\) where \(\mu\) is as above. Following an analogous analysis as in Section~\ref{intro}, there exists a constant \(C_1 > 0\) such that
\begin{equation}\label{an_prob_p}
\mathbb{P}(\hat{A}_n^c(n)) \leq \frac{C_1}{n^2}
\end{equation}
for all \(n \geq 1.\) Finally, set \(\hat{F}_n = \cap_{k=1}^{n} \hat{E}_{k}^c(n) \cap \hat{A}_n(n)\) and fix \(1 \leq k \leq n.\) If \(\hat{F}_n\) occurs, then the time \(\hat{T}^{(k)}_i\) taken to reach \((i,\mathbf{0})\) from \((0,\mathbf{0})\) is less than \(6\mu n,\) for each \(1 \leq i \leq 2n.\) This is because \(t_n(f_i) \geq t_k(f_i)\) and thus \(\hat{A}_n(n) \subset \hat{A}_n(k).\) Since \(\hat{E}^c_{k}(n)\) also occurs, every path \(\pi\) starting from \((0,\mathbf{0})\) and containing \(r \geq \frac{8\mu}{\beta_1}n\) edges has passage time \(\hat{T}_k(\pi)\) at least \(\beta_1 r \geq 8\mu n.\) Therefore, if \(\hat{F}_n\) occurs, the shortest time path with passage time \(\hat{T}^{(k)}_i\) from \((0,\mathbf{0})\) to \((i,\mathbf{0})\) is contained in \(B_{8\mu\beta_1^{-1}n} := [-8\mu\beta_1^{-1}n,8\mu\beta_1^{-1}n]^d\) for each \(1 \leq i \leq 2n\) and for each \(1 \leq k \leq n.\)

From (\ref{a_0k_p}) and (\ref{an_prob_p}), we have that
\begin{equation}\label{fn_h_prob}
\mathbb{P}(\hat{F}_n^c) \leq Cne^{-\beta_2 n} + \frac{C_1}{n^2} \leq \frac{C_2}{n^2}
\end{equation}
for some constant \(C_2 > 0.\) Thus \(\mathbb{P}(\liminf_n \hat{F}_n \cap F_n)= 1.\) \\Fix \(\omega\in\liminf_n~ \hat{F}_n \cap F_n\) and \(m \geq 1.\) For every \(1 \leq k \leq 2m,\) define \(\hat{T}^{(m)}_k = \hat{T}^{(m)}_k(\omega)\) to be the shortest time taken for reaching \((k,0)\) from \((0,0),\) as in Section~\ref{intro}. We have the following result.

\begin{Lemma}\label{etn} We have that
\begin{equation}\label{etn_est}
\mathbb{E}(\hat{T}^{(n)}_n - \mathbb{E}\hat{T}^{(n)}_n)^2 \leq C_1 n^{1+3\alpha}
\end{equation}
for all \(n \geq 1\) and some constant \(C_1 > 0.\)
\end{Lemma}
We prove the above lemma at the end of this section. We use Lemma~\ref{etn} to obtain \(L^2\) convergence of \(\frac{1}{n}\left(T_n - \mathbb{E}T_n\right),\) where \(T_n = T(0,n).\) 
\begin{Corollary} \label{cor2}
\begin{equation}\label{etn_est}
\mathbb{E}(T_n - \mathbb{E}T_n)^2 \leq C_2 n^{\frac{3}{2}-\beta}
\end{equation}
for all \(n \geq 1\) and some positive constants \(C_2\) and \(\beta.\)
\end{Corollary}
\emph{Proof of Corollary~\ref{cor2}}: We have that
\begin{eqnarray}
\mathbb{E}(T_n - \mathbb{E}T_n)^2 \leq 2I_1 + 2\mathbb{E}(\hat{T}^{(n)}_n - \mathbb{E}\hat{T}^{(n)}_n)^2,
\end{eqnarray}
where
\begin{eqnarray}
I_1 &=& \mathbb{E}(T_n - \hat{T}^{(n)}_n - \mathbb{E}(T_n - \hat{T}^{(n)}_n))^2 \nonumber\\
&\leq& 2\mathbb{E}(T_n - \hat{T}^{(n)}_n)^2 + 2 (\mathbb{E}T_n - \mathbb{E}\hat{T}^{(n)}_n)^2 \nonumber\\
&\leq& 4 \mathbb{E}(T_n - \hat{T}^{(n)}_n)^2. \label{i_1}
\end{eqnarray}
It suffices to estimate the last term.

We let \(G_n\) denote the event that the passage time \(t(e_i)\) of every edge in \(B_{8\mu\beta_1^{-1} n}\) is less than \(n^{\alpha}.\) We have that
\begin{equation} \label{t_hat_t_n}
\hat{T}^{(n)}_n \ind(H_n) = T_n \ind(H_n).
\end{equation}
where \(H_n = G_n \cap F_n \cap \hat{F}_n.\)
Thus
\begin{equation} \label{dif_eq}
\mathbb{E}(T_n - \hat{T}^{(n)}_n)^2 = \mathbb{E}(T_n - \hat{T}^{(n)}_n)^2\ind(H_n^c) \leq \left(\mathbb{E}(T_n-\hat{T}^{(n)}_n)^4\right)^{1/2} \left(\mathbb{P}(H_n^c)\right)^{1/2},
\end{equation}
by Cauchy-Schwarz inequality. We have that \[\mathbb{E}(T_n-\hat{T}^{(n)}_n)^4\leq 16\mathbb{E}T_n^4 + 16\mathbb{E}\left(\hat{T}^{(n)}_n\right)^4.\] Since \(T_n \leq \sum_{i=1}^{n} t(f_i),\) where as before, \(f_i\) denotes the edge between \((i-1,0)\) and \((i,0),\)  we have that \[\mathbb{E}T_n^4 \leq n^3\sum_{i=1}^{n}\mathbb{E}t(f_i)^4 \leq C_1n^4\] for some constant \(C_1 > 0.\) An analogous estimate holds for \(\mathbb{E}(\hat{T}_n^{(n)})^4.\) Thus from (\ref{dif_eq}), we have that
\begin{equation}\label{f_n_int}
\mathbb{E}(T_n - \hat{T}^{(n)}_n)^2 \leq C_2n^2 \left(\mathbb{P}(H_n^c)\right)^{1/2},
\end{equation} for some constant \(C_2 > 0.\)

Finally, we choose \(\alpha < \frac{1}{6}\) and \(6(1+d) < K < 6(1+d) + \eta \) such that \(K\alpha > 1+d.\) Here \(\eta > 0\) is as in (iii). We then have that
\begin{equation}\label{g_n}
\mathbb{P}(G_n^c) \leq \sum_{i=1}^{C_3n^d} \mathbb{P}(t(e_i) \geq n^{\alpha}) \leq \frac{C_3n^d}{n^{K\alpha}}\mathbb{E}t(e_i)^{K} \leq \frac{C_4}{n^{1+2\delta}}
\end{equation}
for some positive constants \(C_3,C_4\) and \(\delta.\) Thus from (\ref{fn_prob}), (\ref{fn_h_prob}) and (\ref{f_n_int}), we get that \[\mathbb{E}(T_n - \hat{T}^{(n)}_n)^2 \leq C_5n^2 n^{-\frac{1}{2} - \delta} = C_5n^{\frac{3}{2}-\delta},\] for some positive constant \(C_5.\) \(\qed\)

\emph{Proof of Theorem~\ref{thm1}}: We claim that it suffices to prove that \(\frac{1}{n}\left(\hat{T}^{(n)}_n - \mathbb{E}\hat{T}^{(n)}_n\right)\) converges to zero a.s. Indeed, letting \(H_n \) be as in proof of Corollary~\ref{cor2} and using (\ref{t_hat_t_n}), we have that \[\frac{1}{n}\left(T_n - \mathbb{E}T_n\right) = \frac{1}{n}\left(\hat{T}^{(n)}_n - \mathbb{E}\hat{T}^{(n)}_n\right) + J_{1,n}-\mathbb{E}J_{1,n} - J_{2,n} + \mathbb{E}J_{2,n},\] where \(J_{1,n} = \frac{T_n}{n} \ind(H_n^c)\) and \(J_{2,n} = \frac{\hat{T}^{(n)}_n}{n}\ind(H_n^c).\) From (\ref{fn_prob}), (\ref{g_n}) and Borel-Cantelli Lemma we have that \(\mathbb{P}(\liminf_n H_n) = 1.\) Thus a.s. we have that \(\limsup_n J_{1,n} = 0 = \limsup_n J_{2,n}.\)

It remains to show that \(\mathbb{E}J_{i,n} \rightarrow 0\) as \(n \rightarrow \infty\) for \(i =1,2.\) We show that \(\sup_n\mathbb{E}J_{i,n}^2 < \infty\) for \(i =1, 2.\) This implies that \(J_{1,n}\) and \(J_{2,n}\) are uniformly integrable and completes the claim. We have that \[J_{1,n} \leq \frac{T_n}{n} \leq \frac{1}{n}\sum_{i=1}^{n} t(f_i)\] where as before \(f_i\) denotes the edge from \((i-1,0)\) to \((i,0).\) Thus \[\mathbb{E}J_{1,n}^2 \leq \frac{1}{n}\sum_{i=1}^{n}\mathbb{E}t(f_i)^2 \leq C_1\] for some constant \(C_1 >0\) by condition (iii). An analogous estimate holds for~\(J_{2,n}.\)

To prove that \(\frac{1}{n}\left(\hat{T}^{(n)}_n - \mathbb{E}\hat{T}^{(n)}_n\right)\) converges to zero a.s., we use a subsequence argument as follows. Set \(S_n = \hat{T}^{(n)}_n - \mathbb{E}\hat{T}^{(n)}_n.\) From Lemma~\ref{etn}, we have that \(\mathbb{E}S_n^2 \leq C_1 n^{1+3\alpha}.\) Thus for a fixed \(\epsilon  > 0,\) we have that \[\mathbb{P}(|S_{n^2}| > n^2 \epsilon) \leq \frac{\mathbb{E}S^2_{n^2}}{\epsilon^2 n^4} \leq \frac{C_2}{n^{2-6\alpha}}\] for some constant \(C_2 > 0.\) Since \(\alpha < \frac{1}{6},\) we have that \(2-6\alpha > 1\) and by Borel-Cantelli Lemma, we have that \(\frac{S_{n^2}}{n^2} \rightarrow 0\) a.s. as \(n \rightarrow \infty.\)

We now set \(D_{n^2} = \max_{n^2 \leq k < (n+1)^2} |S_k-S_{n^2}|\) and estimate \(D_{n^2}\) as follows. For \(n^2 \leq k < (n+1)^2,\) we write
\begin{eqnarray}
|S_k - S_{n^2}| &\leq& |\hat{T}^{(k)}_k- \hat{T}^{(n^2)}_{n^2}| + \mathbb{E}|\hat{T}^{(k)}_k - \hat{T}^{(n^2)}_{n^2}| \nonumber\\
&\leq& |\hat{T}^{(k)}_k- \hat{T}^{(k)}_{n^2}|  + |\hat{T}^{(k)}_{n^2}- \hat{T}^{(n^2)}_{n^2}| \nonumber\\
&&\;\;\;\;\;\;\;\;\;\;\; +\;\mathbb{E}|\hat{T}^{(k)}_k- \hat{T}^{(k)}_{n^2}|  + \mathbb{E}|\hat{T}^{(k)}_{n^2}- \hat{T}^{(n^2)}_{n^2}|.\;\;\;\;\label{eq_s2}
\end{eqnarray}
For any integers \(k_1 < k_2 < k_3,\) we have that
\begin{equation}\label{sub_add}
\hat{T}^{(k)}_{k_1,k_3} \leq \hat{T}^{(k)}_{k_1, k_2} + \hat{T}^{(k)}_{k_2,k_3} \text{ and } \hat{T}^{(k)}_{k_1,k_2} \leq \hat{T}^{(k)}_{k_1, k_3} + \hat{T}^{(k)}_{k_2,k_3}.
\end{equation}
Here \(\hat{T}^{(k)}_{k_1,k_2}\) denotes minimum passage time to go from \((k_1,0)\) to \((k_2,0)\) and is defined analogously as \(\hat{T}^{(k)}_n\) for each \(k_1\) and \(k_2.\) Thus \[|\hat{T}^{(k)}_{k} - \hat{T}^{(k)}_{n^2}| \leq \hat{T}^{(k)}_{k,n^2} \leq k^{\alpha}(k-n^2) \leq (n+1)^{2\alpha}((n+1)^2 - n^2) \leq C_1 n^{1+2\alpha}\]  for some constant \(C_1  >0.\) The second inequality is true since the passage time of every edge is less than \(k^{\alpha}.\) Substituting the above estimate into (\ref{eq_s2}), we obtain that
\begin{eqnarray}
|S_k - S_{n^2}| \leq 2C_1 n^{1+2\alpha}  + |\hat{T}^{(k)}_{n^2}- \hat{T}^{(n^2)}_{n^2}| + \mathbb{E}|\hat{T}^{(k)}_{n^2}- \hat{T}^{(n^2)}_{n^2}|.\;\;\;\;\label{eq_s3}
\end{eqnarray}

To estimate the remaining terms, we note that \[0 \leq \hat{T}^{(k)}_{n^2} - \hat{T}^{(n^2)}_{n^2} \leq \hat{T}^{((n+1)^2)}_{n^2} - \hat{T}^{(n^2)}_{n^2}  =: I_{n^2}\] since \(n^2 \leq k < (n+1)^2.\) Thus \[\frac{D_{n^2}}{n^2} \leq \frac{2C_1}{n^{1-2\alpha}} + \frac{I_{n^2}}{n^2} + \frac{\mathbb{E}I_{n^2}}{n^2}.\] We claim that \(\frac{I_{n^2}}{n^2} \rightarrow 0\) a.s. and that \(\frac{I_{n^2}}{n^2}\) is uniformly integrable. Assuming the claims for the moment, we get that  \(\frac{D_{n^2}}{n^2} \longrightarrow 0\) a.s. as \(n \rightarrow \infty.\) For \(n^2 \leq k < (n+1)^2, \) we have that \[\frac{|S_k|}{k}  \leq \frac{|S_k-S_{n^2}|}{k} + \frac{|S_{n^2}|}{k} \leq \frac{|S_k-S_{n^2}|}{n^2} + \frac{|S_{n^2}|}{n^2} \leq \frac{D_{n^2}}{n^2} + \frac{|S_{n^2}|}{n^2}.\] This proves that the original sequence \(\frac{S_k}{k} \rightarrow 0\) a.s. as \(k \rightarrow \infty.\)

To prove the two claims regarding \(I_{n^2},\) we note that \[\hat{T}^{((n+1)^2)}_{n^2}\ind(\hat{H}_n) = \hat{T}^{(n^2)}_{n^2}\ind(\hat{H}_n)\] where \(\hat{H}_n = \hat{F}_{n^2} \cap \hat{F}_{(n+1)^2} \cap \hat{G}_{n^2}\) and \(\hat{G}_{n^2}\) is the event that the passage time \(t(e_i)\) of every edge in \(B_{20\mu\beta_1^{-1} n^2}\) is less than \(n^{2\alpha}.\) As in (\ref{g_n}) we have that \(\mathbb{P}(\hat{G}_{n^2}^c) \leq \frac{C_1}{n^{2+\delta_2}}\) for some constant \(\delta_2  >0.\) From (\ref{fn_h_prob}) and Borel-Cantelli lemma, we then have that \(\mathbb{P}(\liminf_n \hat{H}_n) = 1.\) Since \(I_{n^2} = I_{n^2}\ind(\hat{H}_n^c),\) we get that \(\frac{I_{n^2}}{n^2} \rightarrow 0\) a.s. as \(n \rightarrow \infty.\)

To prove the uniform integrability of \(\frac{I_{n^2}}{n^2},\) we note that \[0 \leq \frac{I_{n^2}}{n^2} \leq \frac{\hat{T}^{((n+1)^2)}}{n^2} \leq \frac{1}{n^2}\sum_{i=1}^{n^2}t(f_i) =: M_n\] where as before \(f_i\) denotes the edge from \((i-1,0)\) to \((i,0).\) Since \(\mathbb{E}M_n^2 \leq \frac{1}{n^2}\sum_{i=1}^{n^2}\mathbb{E}t(f_i)^2 \leq C_1\) for some constant \(C_1 > 0,\) we are done. \(\qed\)

\emph{Proof of Corollary~\ref{cor1}}: We show that \(\frac{\mathbb{E}T(0,n)}{n} \rightarrow \mu\) for some constant \(\mu  >0.\) Since \[\mathbb{E}T(0,n+m) \leq \mathbb{E}T(0,n) + \mathbb{E}T(m,m+n) = \mathbb{E}T(0,n) + \mathbb{E}T(0,m),\] we have by Fekete's Lemma that \[\lim_n\frac{\mathbb{E}T(0,n)}{n} = \inf_{n \geq 1} \frac{\mathbb{E}T(0,n)}{n} =: \mu.\] To show that \(\mu  >0,\) we note that if \(A_{0,k_1}^c\) occurs for \(k_1 = \beta_1 (8\mu)^{-1} n,\) then every path containing \(r \geq 8\mu\beta_1^{-1} k_1 \geq n\) edges has passage time at least \(\beta_1r \geq 8\mu k_1 \geq \beta_1 n.\) Thus \[\mathbb{E}T(0,n) \geq \beta_1 n \mathbb{P}(A^c_{0,k_1}) \geq \beta_3n\] for all \(n \geq 1\) and some constant \(\beta_3 > 0,\) by (\ref{a_0k}). \(\qed\)

\emph{Proof of Lemma~\ref{etn}}: We order the edges as \(e_1,e_2,..\) and for each \( i \geq 1,\) set \({\cal F}_i = \sigma(\hat{t}(e_l) : 1 \leq l \leq i).\) For \(l \geq 1,\) let \(X_l = \mathbb{E}(\hat{T}^{(n)}_n|{\cal F}_l) - \mathbb{E} (\hat{T}^{(n)}_n|{\cal F}_{l-1}).\)  We have that \(0 \leq \hat{T}^{(n)}_n \leq \sum_{i=1}^{n} t_n(f_i) \leq n^{1+\alpha}\) a.s., where as before \(f_i\) denotes the edge from \((i-1,\mathbf{0})\) to \((i,\mathbf{0}).\) Thus we have by Levy's martingale convergence theorem that \[Y_m := \sum_{l=1}^{m} X_l = \mathbb{E}(\hat{T}^{(n)}_n | {\cal F}_{m}) - \mathbb{E}\hat{T}^{(n)}_n \longrightarrow \hat{T}^{(n)}_n - \mathbb{E}(\hat{T}^{(n)}_n)\;\;a.s.\] as \(m \rightarrow \infty.\) By Dominated convergence theorem, we then have that
\[\mathbb{E}(\hat{T}^{(n)}_n -\mathbb{E}\hat{T}^{(n)}_{n})^2 = \mathbb{E}\left(\lim_m Y_m\right)^2 = \lim_m \mathbb{E}Y_m^2.\] By the martingale property, we have that \(\mathbb{E}Y_m^2 = \sum_{l=1}^{m} \mathbb{E}X_l^2.\) We claim that
\begin{equation}
X_l^2 \leq 2n^{2\alpha} \left(\mathbb{P}(e_l \in {\pi}_n | {\cal F}_l) + \mathbb{P}(e_l \in {\pi}_n | {\cal F}_{l-1})\right)\;\;a.s.\label{eq1}
\end{equation}
where \(\pi_n\) is the shortest time path from \((0,0)\) to \((n,0).\) We prove the above result at the end. Using (\ref{eq1}), we obtain that
\begin{eqnarray}
\mathbb{E}(\hat{T}^{(n)}_n -\mathbb{E}\hat{T}^{(n)}_{n})^2 &\leq& 2n^{2\alpha} \sum_{l=1}^{\infty} \mathbb{E}\left(\mathbb{P}(e_l \in {\pi}_n | {\cal F}_{l})  + \mathbb{P}(e_l \in {\pi}_n | {\cal F}_{l-1})\right) \nonumber\\
&=& 4n^{2\alpha}\sum_{l=1}^{\infty} \mathbb{P}(e_l \in {\pi}_n)\nonumber\\
&=& 4n^{2\alpha}\mathbb{E}\sum_{l=1}^{\infty} \ind(e_l \in {\pi}_n)\nonumber\\
&=& 4n^{2\alpha}\mathbb{E}(\#{\pi}_n), \nonumber
\end{eqnarray}
where \(\ind(.)\) refers to the indicator function.

To estimate the length of \(\pi_n,\) let \(\mu = \sup_i \mathbb{E}t(e_i)\) be as in Section~\ref{intro}. We note that if \(\hat{E}^{c}_{k}(n) \) occurs (see paragraph prior to (\ref{a_0k_p})) for \(k \geq \mu^{-1}n^{1+\alpha},\) then every path \(\pi\) with length \(r \geq \frac{8\mu}{\beta_1}k\) has passage time \(\hat{T}_n(\pi)\) at least \(\beta_1 r \geq 8 n^{1+\alpha}.\) Since \(\pi_n\) has passage time at most \(n^{1+\alpha},\) we obtain for \(k \geq \mu^{-1}n^{1+\alpha}\) that \[\mathbb{P}(\#{\pi}_n \geq 8\mu\beta_1^{-1}k) \leq \mathbb{P}(\hat{E}_{k}(n)) \leq e^{-\beta_2 k},\] where \(\beta_2 > 0\) is as in (\ref{a_0k_p}). Since \(\mathbb{E}(\#{\pi}_n) \leq \sum_{k \geq 1} \mathbb{P}(\#{\pi}_n \geq k),\) we obtain that \(\mathbb{E}(\#{\pi}_n) \leq C_1 n^{1+\alpha}\) for some constant \(C_1 > 0.\)

To estimate \(X_l,\) we use the notation of Kesten (1993); for \(j \geq 1,\) let \({\nu}_j(.)\) denote the probability measure associated with \((\hat{t}(e_j),\hat{t}(e_{j+1}),...).\) Let \((\sigma_1,\sigma_2,...)\) and \((\omega_1,\omega_2,...)\) be independent realizations of \((\hat{t}(e_1),\hat{t}(e_2),...)\) and for \(l \geq 1,\) define \([\omega,\sigma]_l = (\omega_1,\omega_2,...,\omega_{l},\sigma_{l+1},\sigma_{l+2},...).\) We have that \[X_l = \int \nu_l(d\sigma) (T_n([\omega,\sigma]_l) - T_n([\omega,\sigma]_{l-1})).\] We note that changing the passage time of edge \(e_l\) does not change the value of the minimum passage time by more than \(n^{\alpha}.\) Also, a change occurs only if \(e_l \in {\pi}_n([\omega,\sigma]_l)\) or \(e_l \in {\pi}_n([\omega,\sigma]_{l-1}).\) Moreover, if Thus \[|\hat{T}^{(n)}_n([\omega,\sigma]_l) - \hat{T}^{(n)}_n([\omega,\sigma]_{l-1})| \leq n^{\alpha} \left(\ind(e_l \in {\pi}_n([\omega,\sigma]_l)) + \ind(e_l \in {\pi}_n([\omega,\sigma]_{l-1}))\right)\] and by Cauchy-Schwarz inequality, we have a.s. that
\begin{eqnarray}
X_l^2 &\leq& \int {\nu}_l(d\sigma)|\hat{T}^{(n)}_n([\omega,\sigma]_l) - \hat{T}^{(n)}_n([\omega,\sigma]_{l-1})|^2 \nonumber\\
&\leq& 2n^{2\alpha} \int {\nu}_l(d\sigma)\left(\ind(e_l \in {\pi}_n([\omega,\sigma]_l)) + \ind(e_l \in {\pi}_n([\omega,\sigma]_{l-1}))\right) \nonumber\\
&=& 2n^{2\alpha} \left(\mathbb{P}(e_l \in {\pi}_n | {\cal F}_l) + \mathbb{P}(e_l \in {\pi}_n | {\cal F}_{l-1})\right). \nonumber
\end{eqnarray}
This proves (\ref{eq1}).\(\qed\)

\section*{Acknowledgements}
I thank Professors Rahul Roy and Thomas Mountford for crucial comments and for my fellowship.



\setcounter{equation}{0} \setcounter{Lemma}{0} \renewcommand{\theLemma}{II.%
\arabic{Lemma}} \renewcommand{\theequation}{II.\arabic{equation}} %
\setlength{\parindent}{0pt}




%





\bibliographystyle{plain}

\end{document}